\documentclass[12pt, final]{amsproc}
\usepackage[cp850]{inputenc}
\usepackage{amssymb}
\usepackage[mathscr]{eucal}
\usepackage{amscd}
\usepackage{amsmath}
\usepackage{setspace}
\usepackage{tikz}
\usetikzlibrary{shapes.geometric}
\def\Ker{\operatorname{Ker}}

\def\dim{\operatorname{dim}}
\def\supp{\operatorname{supp}}

\theoremstyle{plain}

\newtheorem{theorem}{Theorem}
\newtheorem{proposition}{Proposition}
\newtheorem{lemma}{Lemma}
\newtheorem{corollary}{Corollary}

\newtheorem*{theorem2}{Kalton uniqueness theorem}
\theoremstyle{remark}
\newtheorem{remark}{Remark}

\title{SOME MORE TWISTED HILBERT SPACES}

\address{Departamento de Matem\'aticas, Universidad de Extremadura, Avenida de Elvas s/n, 06011 Badajoz, Spain.}
\email{danmorg@unex.es}
\author[D. Morales]{Daniel Morales}
\address{Departamento de Matem\'aticas, Universidad de Extremadura, Avenida de Elvas s/n, 06011 Badajoz, Spain.}
\email{jesus@unex.es}
\author[J. Su\'arez]{Jes\'us Su\'arez}


\thanks{The first author was supported by the grant BES-2017-079901 of the project  MTM2016-76958-C2-1-P. The second author was supported in part by projects MTM2016-76958-C2-1-P, PID2019-103961GB-C21 and IB16056. This is part of the thesis of first named author under the supervision of Jes\'us M.F. Castillo and the second named author.}

\begin{document}

\begin{abstract} We provide three new examples of twisted Hilbert spaces by considering properties that are ``close" to Hilbert. We denote them $Z(\mathcal J)$, $Z(\mathcal S^2)$ and $Z(\mathcal T_s^2)$. The first space is asymptotically Hilbertian but not weak Hilbert. On the opposite side, $Z(\mathcal S^2)$ and $Z(\mathcal T_s^2)$ are not asymptotically Hilbertian. Moreover, the space $Z(\mathcal T_s^2)$ is a HAPpy space and the technique to prove it gives a ``twisted" version of a theorem of Johnson and Szankowski (Ann. of Math. 176:1987--2001, 2012). This is, we can construct a nontrivial twisted Hilbert space such that the isomorphism constant from its $n$-dimensional subspaces to $\ell_2^n$ grows to infinity as slowly as we wish when $n\to \infty$. 
\end{abstract}


\maketitle
\section{Introduction}
 Since its inception as a solution to the Palais problem \cite{ELP}, twisted Hilbert spaces, i.e. Banach spaces $X$ admitting a Hilbert subspace $H$ so that the corresponding quotient $X/H$ is Hilbert, have been a fruitful place where to seek counterexamples. And still now there are natural problems in the literature, as for example the ergodicity problem \cite{FR}, for which nontrivial twisted Hilbert spaces could be a solution. The scarcity of a broad variety of twisted Hilbert spaces has been a problem in the area. Indeed, for a long time, the only known examples were the Enflo-Lindenstrauss-Pisier space (\cite{ELP}) that we denote $\ell_2(\mathcal E_n)$, and the Kalton-Peck space $Z_2$ (\cite{KaPe}). This paper provides a few more examples, somehow continuing the work in \cite{derivation,JS1} and \cite{JS2}.

To put the forthcoming results in perspective, let us recall that twisted Hilbert spaces are ``close" to Hilbert spaces in many senses: they are $\ell_2$-saturated, superreflexive and they have type $2-\varepsilon$ and cotype $2+\varepsilon$ for all $\varepsilon>0$. However, there are boundaries: Maurey's extension theorem shows us that nontrivial twisted Hilbert spaces cannot have type $2$ or cotype $2$ (cf. \cite{Mau}); and a deep result of Kalton shows that they cannot even have an unconditional basis \cite{ka2}. Nevertheless, regarding twisted Hilbert spaces with extremal properties, the second author constructed in \cite{JS1} a twisted Hilbert space $Z(\mathcal T^2)$, where $\mathcal T^2$ denotes the $2$-convexification of the Tsirelson space $\mathcal T$, that is a weak Hilbert space. This last space and the Kalton-Peck space $Z_2$ emerge from the same scheme: Let $X$ be a separable Banach space for which complex interpolation yields a Hilbert space in the form $(X, X^*)_{1/2} = \ell_2$. Then, a twisted Hilbert space $Z(X)$ arises as the derivation of the previous formula. It is known as the derived space of the interpolation space at $1/2$. In particular, we have $Z(\ell_1)=Z_2$ while $Z(\mathcal T^2)$ is weak Hilbert. 

In what follows, we will focus on derived spaces $Z(X)$ sharing properties ``close" to Hilbert, showing that there is still room enough until the Hilbert barrier. Beyond weak Hilbert, notions of this type found in the literature are: asymptotically Hilbertian (as. Hilbertian) spaces, spaces with the property $(H)$, different forms of $E(n,m,K)$-properties as introduced by Nielsen and Tomczak-Jaegermann \cite{NT} and HAPpy spaces, a term coined by Johnson and Szankowski in \cite{JS}.

 In this context, our first example is an asymptotically Hilbertian space $Z(\mathcal J)$ without the property $(H)$. Hence, it is not a weak Hilbert space. For the specialist it is perhaps not surprising that this is achieved with the right choice of $p_n, k_n$ in $\mathcal J=\ell_2(\ell_{p_n}^{k_n})$. Our second example picks as $X$ the $2$-convexification $\mathcal S^2$ of the Schreier space $\mathcal S$, see e.g. \cite{CS}. We show that $Z(\mathcal S^2)$ is not asymptotically Hilbertian but its natural basis has the $E(n,n,K)$-property (see Section \ref{background} for the precise definitions) while the basis of our first and third examples lack it.
 This third example is modelled over $\mathcal T_s^2$ which is the Casazza-Nielsen symmetric version of $\mathcal T^2$, cf. \cite{CN2}. The distance of the $n$-dimensional subspaces of $Z(\mathcal T_s^2)$ to Hilbert grows very slowly to infinity. This is in contrast with $Z(\mathcal S^2)$ that contains an isomorphic copy of $Z_2$ in spite of the $E(n,n,K)$-property; a second turn of the screw will show that $Z(\mathcal S^2)$ is not isomorphic to a subspace of $Z_2$.

As a by-product of our methods, we may give a ``twisted" version of a result of Johnson and Szankowski \cite[Theorem 3.1]{JS}, namely, given $1<\delta_n\to \infty$ there is a twisted Hilbert $Z$ so that $d_n(Z)\leq \delta_n$, where $d_n(Z)$ is the supremum over all the $n$-dimensional subspaces $E$ of $Z$ of the isomorphism constant from $E$ to $\ell_2^n$. This is a way to construct nontrivial examples of HAPpy twisted Hilbert spaces. Previously, the only known example of such kind was $Z(\mathcal T^2)$, but this is weak Hilbert and thus as. Hilbertian. However, the space $Z$ constructed above as well as $Z(\mathcal T_s^2)$ are not as. Hilbertian.

Of course, the spaces $Z(\mathcal J)$, $Z(\mathcal S^2)$ and $Z(\mathcal T_s^2)$ are mutually non isomorphic and are also not isomorphic to the known previous examples of twisted Hilbert spaces: $\ell_2(\mathcal E_n)$, $Z_2$ and $Z(\mathcal T^2)$. In a broader sense, our examples are representatives of three new (and very natural) categories of twisted Hilbert spaces.

The paper is organized as follows. The next section contains a short description of the background needed. Sections \ref{ahnw},\ref{ex2},\ref{cinco} are devoted to describe each example separately. We have included one last section with a picture that might help the reader to organize our examples of twisted Hilbert spaces.
\section{Background}\label{background}
We use standard notation for Banach spaces as provided in the book of Albiac and Kalton \cite{AK}. We reserve the word space for Banach space, either finite or infinite dimensional. In the finite dimensional setting we will very often write $$d_E=d(E,\ell_2^{\dim E}),$$ where $d$ stands for the Banach-Mazur distance. In this sense, given a space $X$, we also define $$d_n(X)=\sup d_E,$$ where the supremum runs over all $n$-dimensional subspaces $E$ of $X$. Recall that for a space $X$ the number $a_{n,2}(X)$ is defined to be the least constant $a$ such that
$$\mathbb E\left \|\sum_{j=1}^n \varepsilon_j x_j\right\|\leq a \left ( \sum_{j=1}^n \| x_j\|^2\right)^{1/2},$$
for all $x_1,...,x_n\in X$ and where the average is taken over all $(\varepsilon_j)_{j=1}^n\in \{-1,1\}^n$. The space $X$ has \textit{type $2$} if $a_2(X):=\sup_{n\in \mathbb N}a_{n,2}(X)<\infty$. The cotype 2 constant $c_{n,2}(X)$ is defined in a similar vein as the least constant $c$ such that $$\left ( \sum_{j=1}^n \| x_j\|^2\right)^{1/2}\leq c\cdot \mathbb E\left \|\sum_{j=1}^n \varepsilon_j x_j\right\|,$$
for all $x_1,...,x_n\in X$ and thus $X$ has \textit{cotype $2$} if $c_2(X):=\sup_{n\in \mathbb N}c_{n,2}(X)<\infty$. As an example, let us recall that $\ell_p$ has type $\min \{p,2\}$ and cotype $\max\{p,2\}$, see \cite{SM}.
A remarkable fact due to Kwapie\'n \cite{Kw} is that
$$d_E\leq a_2(E)\cdot c_2(E),$$
which will be used throughout the paper.
\subsection{Background on Hilbert-like properties}
 We will mainly consider five notions related to a Hilbert space: Weak Hilbert spaces, as. Hilbertian, the property $(H)$, the $E(n,n,K)$-property and the HAPpy spaces.

We say $X$ is a weak Hilbert space if there is $0<\delta_0<1$ and a constant $C$ with the following property: every finite dimensional subspace $E$ of $X$ contains a subspace $F\subseteq E$ with $\dim F\geq \delta_0 \dim E$ such that $d_F\leq C$ and there is a projection $P:X\to F$ with $\|P\|\leq C$.
The definition above is not the original one but is chosen out among the many equivalent characterizations given by Pisier \cite[Theorem 12.2.(iii)]{Pi}. Recall that $\mathcal T^2$ is the most classical example of a nontrivial weak Hilbert space and the reader may find in \cite{CS} or \cite{Pib} a comprehensive study.

Closely related to this is the notion of a space $X$ that is as. Hilbertian, meaning that there is $C>0$ such that for every $n\in \mathbb N$, there is a finite codimensional $X_n$ in $X$ so that every $n$-dimensional subspace $E$ of $X_n$ satisfies $d_E\leq C$. 


Just in between, lies the notion of the property $(H)$. We say that $X$ has the property $(H)$  if for each $\lambda>1$ there is a $K(\lambda)$ such that for any $n\in \mathbb N$ and any $\lambda$-unconditional normalized (basic) sequence $(u_j)_{j=1}^n\in X$, we have
$$K(\lambda)^{-1}\sqrt{n}\leq \left \| \sum_{j=1}^n  u_j \right\|\leq K(\lambda)\sqrt{n}.$$
Using the characterizations of weak Hilbert spaces given by Pisier, it is easy to prove that they have the property $(H)$, see \cite[Proposition 14.2.]{Pib}. A bit more elaborated is the proof that the property $(H)$ implies as. Hilbertian, see \cite[Theorem 14.4]{Pib}.

We have left almost at the end the $E(n,n,K)$-property. While the idea behind this notion is simple, it requires two twists. Let $X$ be a Banach space with a basis $(e_j)_{j=1}^{\infty}$ and let $K\geq 1$ and $m,n\in \mathbb N$. First, we say  that $(e_j)_{j=1}^{\infty}$ is $(m,K)$-euclidean, if for every subset $A\subseteq \mathbb N$ with $|A|\leq m$, $(e_j)_{j\in A}$ is $K$-equivalent to the
unit vector basis of $\ell_2^{|A|}$. Then $(e_j)_{j=1}^{\infty}$ is said to have the $E(n, m,K)$-property if there is a set $I\subseteq \mathbb N$, with $|I|=n$, so that $$\{e_j: j\in \mathbb N \backslash I  \}$$ is $(m,K)$-euclidean. A key result for us is that for any unconditional sequence $(u_j)_{j=1}^{\infty}$ with the property $(H)$ there is $K>0$ such that $(u_j)_{j=1}^{\infty}$ satisfies the $E(n,n,K)$-property for every $n\in \mathbb N$. This can be found in the paper of Nielsen and Tomczak-Jaegermann \cite[Proposition 3.8]{NT}. A particular case of the $E(n,m,K)$-property is the $E_t(n,m,K)$-property, meaning that $\{ e_j : j\geq n+1\}$ is $(m,K)$-euclidean.

Finally, following Johnson and Szankowski \cite{JS}, we say that a space $X$ has the hereditary approximation property (HAP) or is a HAPpy space if all the subspaces of $X$ have the approximation property.
\subsection{Background on complex interpolation and twisted Hilbert spaces}
Let $\omega$ denote the vector space of complex scalar sequences endowed with the pointwise convergence. Let $X$ be a space with an unconditional basis. An homogeneous map $\Omega:X\longrightarrow \omega$ is called a \textit{centralizer} if there is a constant $C$ so that:
\begin{equation}\label{centralizer}
\| \Omega(ax)-a\Omega(x)\|_X\leq C\|a\|_{\infty}\|x\|_X,\;\;\:\;a\in \ell_{\infty},x\in X.
\end{equation}
A typical way to obtain centralizers is through complex interpolation \cite{Cald} and the classic reference for interpolation is the book of Bergh and L\"ofstr\"om \cite{BeLo}. We only describe some basic facts. Let $X_0,X_1$ be a couple of spaces with a joint $1$-unconditional basis and natural inclusions into $\omega$. We shall consider the vector space $\mathcal F_{\infty}(X_0,X_1)$ of all functions $F:\mathbb S\to \omega$, which are bounded and continuous on the strip $$\mathbb S=\{ z: 0\leq Rez \leq 1 \},$$
and analytic on the open strip $\mathbb S_0=\{z: 0< Re z <1 \},$ and moreover, the functions $t\in \mathbb R\to F(j+it)\in X_j$ with $j=0,1$ are bounded and continuous functions. The vector space $\mathcal F_{\infty}(X_0,X_1)$ is a Banach space when is endowed with the norm $$\|F\|_{\mathcal F_{\infty}}=\max \left ( \sup_{t\in \mathbb R} \| F(it)\|_{X_0},   \sup_{t\in \mathbb R} \| F(1+it)\|_{X_1}\right).$$
 The interpolation space $X_{\theta}:=(X_0,X_1)_{\theta}$ consists of all $x\in \omega$ such that $x=F(\theta)$ for some $F\in \mathcal F_{\infty}(X_0,X_1)$ endowed with the quotient norm $$\|x\|_{X_{\theta}}=\inf \{ \|F\|_{\mathcal F_{\infty}(X_0,X_1)} : F(\theta)=x\}.$$
  We denote as usual $\delta_{\theta}: \mathcal F_{\infty}(X_0,X_1) \to  X_{\theta}$  the natural quotient map where $\theta\in(0,1)$.
 Fix now a constant $\rho>1$ and thus for each $x\in X_{\theta}$ pick a map $B(x)\in \mathcal F_{\infty}(X_0,X_1)$ with $B(x)(\theta)=x$ and $\|B(x)\|_{\mathcal F_{\infty}}\leq \rho\|x\|_{\theta}$ for which there is no loss of generality in assuming that it is homogeneous. A centralizer $\Omega:X_{\theta}\longrightarrow \omega$ comes defined as $$\Omega(x)=\delta'_{\theta}B(x).$$
Let us explain the connection between centralizers and twisted Hilbert spaces. Recall that a \textit{short exact sequence} is a diagram like
\begin{equation*}\label{ses}
\begin{CD} 0 @>>>Y@>j>>Z@>q>>X@>>>0
\end{CD}
\end{equation*} where the morphisms are linear and continuous and such that the image of each arrow is the kernel of the next one. This condition implies that $Y$ is a subspace of $Z$ through $j$ and thanks to
the open mapping theorem we find that $X$ is isomorphic to $Z/j(Y)$. We usually refer to $Z$ as a \textit{twisted sum} of $Y$ and $X$ and if $Y\approx\ell_2\approx X$ we simply say that $Z$ is a \textit{twisted Hilbert space}.
The interpolation scheme described before produces a natural twisted Hilbert $Z(X)$ of $(X, X^*)_{1/2}=\ell_2$ called \textit{the derived space} defined as the set of couples $(x,y)\in \omega\times \omega$ for which the following norm $$\|(x,y)\|_{Z(X)}=\inf \{ \|F\|: F(\theta)=y, F'(\theta)=x\}$$ is finite.
Just as an explanation, the derived space is denoted by $d(X,X^*)_{1/2}$ in the survey of Kalton and Montgomery-Smith \cite{KM}. It will take an instant to the reader to check (see the discussion in \cite[Page 1159]{KM}) that we have a short exact sequence
$$
\begin{CD} 0 @>>>\ell_2@>j>>Z(X)@>q>>\ell_2@>>>0,
\end{CD}
$$
where $j(x)=(x,0)$ and $q(x,y)=y$. This twisted Hilbert admits a representation in terms of a centralizer $\Omega$. The exact relationship is that the norm of $Z(X)$ is equivalent to the quasi-norm 
\begin{equation}\label{quasinorma}
\|(x,y)\|=\|x-\Omega(y)\|+\|y\|,
\end{equation} where $\Omega$ is, of course, a centralizer corresponding to  $(X, X^*)_{1/2}=\ell_2$, see \cite[Proposition 7.2.]{CCS}. The perfect example where all these ideas crystallize is the Kalton-Peck space which corresponds to $Z(\ell_1)=Z_2$. It is well known that the centralizer comes defined for a norm one vector $x=\sum_{j=1}^{\infty}x_je_j$ as 
\begin{equation}\label{kpmap}
\Omega(x)=\sum_{j=1}^{\infty}x_j\log|x_j|e_j,
\end{equation} with the agreement that $0\cdot \log 0=0$. The map (\ref{kpmap}) is the so-called Kalton-Peck map \cite{KaPe} for which it is helpful to write $\Omega(x)=x\cdot \log x$ with the obvious meaning.

Let us recall now a couple of results that will be useful for us. The first is that $Z(X)$ is isomorphic to its own dual \cite[Lemma 1]{JS2} while the second is the existence of a natural basis in $Z(X)$.
\begin{proposition}\label{lambda} Let $(e_j)_{j=1}^{\infty}$ be the canonical unconditional basis of $(X, X^*)_{1/2}=\ell_2$ and set $v_{2j-1}=(e_j,0), v_{2j}=(0,e_j)$ for $j\in \mathbb N$. Then
\begin{enumerate}
\item[$(i)$] The sequence $(v_j)_{j=1}^{\infty}$ is a basis for $Z(X)$.
\item[$(ii)$] The sequence $(v_{2j})_{j=1}^{\infty}$ is unconditional.
\end{enumerate}
\end{proposition}
The first part is proved by adapting the proof of \cite[Theorem 4.10]{KaPe} while the second follows picking $a\in \{-1,1\}^{n}$ in (\ref{centralizer}).
At this stage, the reader may wonder whether $Z(X)$ is again Hilbert or not. The answer to this fair question is given by the so-called Kalton uniqueness theorem \cite[Theorem 7.6]{ka}.
\begin{theorem2}
Let $(X,X^*)$ and $(Y_0,Y_1)$ be two couples of spaces all of them with a joint unconditional basis $(e_j)_{j=1}^{\infty}$. Assume that $(X,X^*)_{1/2}=\ell_2$ with corresponding centralizer $\Omega_X$ and  $(Y_0,Y_1)_{1/2}=\ell_2$ with centralizer $\Omega_Y$. If $\Omega_X$ and $\Omega_Y$ are equivalent, then $X$ and $Y_0$ have equivalent norms.
\end{theorem2}
The interested reader will find out that  \cite[Theorem 7.6]{ka} has two parts: existence and uniqueness. While the existence part contains a technical assumption on the centralizer under discussion, the uniqueness part does not. This can be easily checked in the first five lines of proof. Indeed, the proof only observes that the claim on the centralizers is a claim on the indicator functions involved and then we may apply \cite[Proposition 4.5]{ka}.

Kalton is claiming that the couples are determined by the corresponding derivation, up to the natural equivalent relation of centralizers which is that the difference is bounded. This is, two centralizers $\Omega_X$ and $\Omega_Y$ are \textit{equivalent} if the difference $\Omega_X-\Omega_Y$ is an (homogeneous) bounded map $\ell_2\to \omega$. In particular, a centralizer is \textit{bounded} if it is equivalent to zero. A direct consequence of this in our setting is the formal answer to the previous question.
\begin{corollary} Let $(X,X^*)$ be spaces with a joint unconditional basis so that $(X,X^*)_{1/2}=\ell_2$. The following conditions are equivalent.
\begin{enumerate}
\item The spaces $X$ and $\ell_2$ have equivalent norms.
\item The twisted-Hilbert $Z(X)$ is isomorphic to $\ell_2$.
\item The space $[v_{2j}]_{j=1}^{\infty}$ is isomorphic to $\ell_2$.
\end{enumerate}
\end{corollary}
\begin{proof}
(1) $\Rightarrow$ (2) If the norm of $X$ is equivalent to $\ell_2$ then so does the norm of $X^*$. Then, for each $x\in \ell_2$, the constant function $B(x)(z)=x$ is a bounded selection (where the bound depends only on the previous equivalence constants) whose derivative is zero. Thus, since $\Omega_X=0$, a simple appeal to the formula (\ref{quasinorma}) finishes.

(2) $\Rightarrow$ (3) Obvious.

(3) $\Rightarrow$ (1) If $[v_{2j}]_{j=1}^{\infty}$ is isomorphic to $\ell_2$, then $[v_{2j}]_{j=1}^{\infty}$ must have type $2$. In particular, we have that there is $C>0$ so that $$\mathbb E \left \|\left(0, \sum_{j=1}^{\infty}\varepsilon_ja_je_j \right)\right\|\leq C\left(\sum_{j=1}^{\infty}a_j^2 \right)^{1/2},$$
for every $(a_j)_{j=1}^{\infty}\in \ell_2$. Since the basis $\{(0,e_j)\}_{j=1}^{\infty}$ is $K$-unconditional by Proposition \ref{lambda}(ii), we readily find from the estimate above that
$$\left \|\left(0, \sum_{j=1}^{\infty}a_je_j \right)\right\|\leq KC\left(\sum_{j=1}^{\infty}a_j^2 \right)^{1/2},\;\;(a_j)_{j=1}^{\infty}\in \ell_2.$$
Using the expression (\ref{quasinorma}), we deduce from above that $$\left \|\Omega_X\left (\sum_{j=1}^{\infty}a_je_j \right)\right\|\leq (KC-1)\left(\sum_{j=1}^{\infty}a_j^2 \right)^{1/2},\;\;(a_j)_{j=1}^{\infty}\in \ell_2,$$
which is to say that $\Omega_X$ is bounded. This means, by definition, that $\Omega_X$ is equivalent to the centralizer $0$ that is represented by the couple $(\ell_2,\ell_2)$ as we argued in (1) $\Rightarrow$ (2). Kalton's result enters now into the game and shows that the norms of $X$ and $\ell_2$ must be equivalent.
\end{proof}
There is a straight and elementary route to this corollary which avoids Kalton's uniqueness theorem. The proof replaces Kalton's result by the method of critical points which is a new and simple way to describe all centralizers. This method may be found in the forthcoming paper \cite{JS3}.
\section{An asymptotically Hilbertian not weak Hilbert space}\label{ahnw}
Let us define for each $n=1,...$ the numbers $k_n=2^{n+1}$ and 
\begin{equation}\label{empezamos}
\frac{1}{p_n}-\frac{1}{2}=\frac{1}{\sqrt{n}}=\frac{1}{2}-\frac{1}{p_n^*}.
\end{equation}
Let $\mathcal J=\ell_2(\ell_{p_n}^{k_n})$ and $\mathcal J^*=\ell_2(\ell_{p_n^*}^{k_n})$, so that $$(\mathcal J,\mathcal J^*)_{1/2}=\ell_2,$$
see for example \cite{BeLo}. Let us recall that $Z(\mathcal J)$ denotes the induced twisted Hilbert space.
\begin{theorem}\label{as} 
The space $Z(\mathcal J)$ is asymptotically Hilbertian but fails the $E(n,n,K)$-property. In particular, it fails the property  $(H)$ and it is not a weak Hilbert space.
\end{theorem}
\begin{proof}
The spaces $\mathcal J$ and $\mathcal J^*$ are non Hilbertian examples of as. Hilbertian. This seems to be folklore but we have been unable to find a proof in the literature so we sketch one for the sake of clarity. 
\begin{flushleft}
\textbf{Claim A:} The space $\mathcal J$ is asymptotically Hilbertian.
\end{flushleft}
To see that $\mathcal J$ is not Hilbert one just need to recall that
\begin{equation*}
d(\ell_{p_n}^{k_n},\ell_2^{k_n})=k_n^{\frac{1}{p_n}-\frac{1}{2}}=2^{\sqrt{n}}\cdot 2^{1/\sqrt{n}}\to \infty,
\end{equation*}
as $n\to \infty$; the same argument shows that $\mathcal J^*$ is not Hilbert. To check that $\mathcal J$ is as. Hilbertian observe that we have, for $m\geq n^2$, the following set of estimates:
\begin{eqnarray*}
a_{k_n,2}(\ell_{p_m}^{k_m})&\leq & k_n^{\frac{1}{p_m}-\frac{1}{2}}\\
&\leq& k_n^{\frac{1}{p_{n^2}}-\frac{1}{2}}\\
&=& 2^{1/n}\cdot 2\\
&\leq& 4,
\end{eqnarray*}
where the first inequality follows, for example, from a simple inspection of the proof of \cite[9.3. Example]{SM} and Kahane's inequality \cite[9.2.]{SM}. Thus, we trivially find
\begin{equation}\label{dificil}
a_{k_n,2}\left( \oplus_{m=n^2}^{\infty} \ell_{p_m}^{k_m} \right)_2\leq 4.
\end{equation}
Let us denote for simplicity $\mathcal J_{n^2}=\left( \oplus_{m=n^2}^{\infty} \ell_{p_m}^{k_m} \right)_2$.
In particular, if $E\subseteq \mathcal J_{n^2}$ with $\dim E\leq k_n$, then
\begin{eqnarray*}
d_E&\leq& a_{2}(E)\cdot c_{2}(E)\\
&\leq& 2\sqrt{2 \pi}\cdot a_{k_n,2}(\mathcal J_{n^2})\cdot c_{k_n,2}(\mathcal J_{n^2})\\
&\leq& 8\sqrt{2 \pi}\cdot c_{k_n,2}(\mathcal J_{n^2})\\
&\leq& 8\sqrt{2 \pi}\cdot c_2,
\end{eqnarray*}
for some absolute $c_2>0$ since $\mathcal J$ (and thus $\mathcal J_{n^2}$) has cotype 2. The first inequality is the remarkable fact due to Kwapie\'n \cite{Kw} of the introduction while the second holds true because of a result of Tomczak-Jaegermann \cite[Theorem 2]{To}. One last remark is that $\mathcal J_{n^2}$ has codimension $\sum_{m=1}^{n^2-1} k_m=2(2^{n^2}-2)$, so the case of $\mathcal J$ is done. A similar argument works for $\mathcal J^*$ interchanging the roles of type and cotype. Indeed, in this last case, we have that
\begin{equation}\label{facil}
a_{k_n,2}\left( \oplus_{m=n^2}^{\infty} \ell_{p_m^*}^{k_m} \right)_2\leq a_2,
\end{equation}
for some absolute $a_2>0$ since $\mathcal J^*$ has type 2.
\begin{flushleft}
\textbf{Claim B:} The twisted-Hilbert $Z(\mathcal J)$ is asymptotically Hilbertian.
\end{flushleft}
We write as usual $Z(\mathcal J_{n^2})$ the induced twisted Hilbert. Then, we may invoke \cite[Proposition 3]{JS1} with the estimates (\ref{dificil}) and (\ref{facil}) to conclude that
\begin{equation}\label{facil2}
a_{k_n,2}\left( Z(\mathcal J_{n^2}) \right)\leq C,
\end{equation}
for some absolute $C>0$. Let us recall that $Z(\mathcal J_{n^2})$ is $\lambda$-isomorphic to its dual, for some constant $\lambda>0$ independent of $n$; a detailed proof may be found in \cite[Lemma 1]{JS2}. Hence, it is easy to deduce that
\begin{equation}\label{facil3}
c_{k_n,2}\left( Z(\mathcal J_{n^2}) \right)\leq \lambda \cdot C,
\end{equation}
by using an argument contained in \cite[Proposition 3.2]{Pi2}. Therefore, if $E\subseteq Z(\mathcal J_{n^2})$ denotes a subspace with $\dim E\leq k_n$, then we have by (\ref{facil2}) and (\ref{facil3})
\begin{eqnarray*}
d_E&\leq& a_2(E)\cdot c_2(E)\\
&\leq&2\sqrt{2\pi}\cdot a_{k_n,2}\left( Z(\mathcal J_{n^2}) \right)\cdot c_{k_n,2}\left( Z(\mathcal J_{n^2}) \right)\\
&\leq& 2\sqrt{2\pi}\cdot \lambda \cdot C^2.
\end{eqnarray*}
One last remark is that $Z(\mathcal J_{n^2})$ is a subspace of $Z(\mathcal J)$ of codimension $2\cdot \sum_{m=1}^{n^2-1} k_m=2^2(2^{n^2}-2)$.
\begin{flushleft}
\textbf{Claim C:} If $Z(\mathcal J)$ satisfies the $E(2^n,2^n,K)$-property, then $2K\geq \sqrt{n}$.
\end{flushleft}
Let us consider $Z(\ell_{p_n}^{k_n})$ for each $n=1,...$, so that $Z(\mathcal J)=\ell_2(Z(\ell_{p_n}^{k_n}))$. Assume $Z(\mathcal J)$ has the $E(n,n,K)$-property so that the unconditional sequence $$U=\{(0,e_j)\in Z(\ell_{p_n}^{k_n}): j=1,...,k_n; n=1,...  \}$$ also has, for some $K>0$, the $E(2^n,2^n,K)$-property. Fix $n$ and pick any set of $2^n$-vectors in $U$, say $A$. Observe that the number of vectors $(0,e_j)$ in each $Z(\ell_{p_m}^{k_m})$ is exactly $k_m=2\cdot 2^m$ for each $m\in \mathbb N$. So, no matter how we choose $A$, we still have a disjoint set with $A$, say $B$, of $2^n$ vectors of the type $(0,e_j)$ in $Z(\ell_{p_n}^{k_n})$. This set of vectors must be, by definition of the $E(2^n,2^n,K)$-property, $K$-equivalent to the basis of $\ell_2^{|B|}$. Recall that the centralizer $\Omega_n$ corresponding to $Z(\ell_{p_n}^{k_n})$ is of the form
 $$\Omega_n=\left(\frac{2}{p_n^*}-\frac{2}{p_n}\right)\Omega:\ell_2^{k_n}\longrightarrow \omega,$$
 where $\Omega(x)=x\cdot \log x$ is the Kalton-Peck map for $\|x\|=1$, cf. \cite[Page 1160]{KM}. Therefore, we find that the following must hold
 \begin{equation}\label{Omega1}
\left\|\Omega_n\left(\sum_{j\in B} e_j\right)\right \|\leq (K-1)\left \| \sum_{j\in B} e_j \right\|=(K-1)\sqrt{2^n}.
\end{equation}
But on the other hand
\begin{equation}\label{Omega2}
\left\|\Omega_n\left(\sum_{j\in B} e_j\right)\right \|=\frac{2^2}{\sqrt{n}}\cdot\sqrt{2^n}\log \sqrt{2^n}=2\cdot \sqrt{n}\cdot \sqrt{2^n}\log2,
\end{equation}
where we have used that
$$\frac{1}{p_n}-\frac{1}{p_n^*}=\frac{1}{p_n}-\frac{1}{2}+\frac{1}{2}-\frac{1}{p_n^*}=\frac{2}{\sqrt{n}},$$
provided by our choice in (\ref{empezamos}). If we plug (\ref{Omega2}) in (\ref{Omega1}) and simplify, we reach to
$$\sqrt{n}\cdot 2\log2\leq K-1,$$
so the claim is proved.
\end{proof}




\section{A non asymptotically Hilbertian space with the $E(n,n,K)$-property}\label{ex2}
 A natural precursor of the Tsirelson space $\mathcal T$ is the Schreier space $\mathcal S$ which is defined as the completion of $c_{00}$, the vector space of finitely supported sequences, under the following norm
$$\|x\|_{\mathcal S}=\sup_{A\in \mathbb A}\sum_{j\in A} |x_j|,$$
where $\mathbb A$ denotes the set of admissible subsets of $\mathbb{N}$; recall that a finite subset $A=\{n_1<...<n_k\}$ is \textit{admissible} if $k\leq n_1$. 
If we denote $\mathcal S^2$ the 2-convexification of $\mathcal S$, then we have that $(\mathcal S^2,(\mathcal S^2)^*)_{\frac{1}{2}}$ and $\ell_2$ have equivalent norms by \cite[Corollary 4.3]{CoS}. Therefore, the derived space $Z(\mathcal S^2)$ is a twisted Hilbert space 
that somehow plays the role of a natural precursor of $Z(\mathcal T^2)$. In this sense, the basis of $Z(\mathcal S^2)$ undergoes a property typically shared by weak Hilbert spaces with an unconditional basis: the $E(n,n,K)$-property \cite[Proposition 3.8.]{NT}.
\begin{proposition}\label{Enn} There is $K>0$ such that the basis $(v_j)_{j=1}^{\infty}$ of $Z(\mathcal S^2)$ has the $E(n,n,K)$-property for all $n\in \mathbb N$.
\end{proposition}
\begin{proof} Let us prove first that `there is $K>0$ such that the subsequence $(v_{2j})_{j=1}^{\infty}$ has the $E_t(n,n,K)$-property for every $n\in \mathbb N$'. Given $y\in c_{00}$ with $\supp(y)\in \mathbb A$, we find
\begin{equation}\label{cdelta2}
 \max(\|y\|_{\mathcal S^2},\|y\|_{(\mathcal S^2)^*})=\|y\|_{\ell_2}. 
\end{equation}
Therefore, given $y\in c_{00}$ as above, let us consider the constant function $$B(y)(z)= y\in \mathcal F_{\infty}(\mathcal S^2,(\mathcal S^2)^*).$$
 We have by (\ref{cdelta2}) that $B(y)$ is a $1$-bounded selection for $y$ (see \cite[p. 1159]{KM}) and thus the centralizer $\Omega(y):=\delta'_{1/2}B(y)$ must be zero.  Therefore, for $y\in c_{00}$ with $\supp (y)\in \mathbb A$, we have that
$$\left \|\Omega(y) \right \| + \|y \|=\left \|y \right \|.$$
The expression above shows that for any $A\in \mathbb A$, we have that $(v_{2j})_{j\in A}$ is equivalent to the basis of $\ell_2^{|A|}$. This equivalence holds in the quasi-norm induced by $\Omega$ that is equivalent to the norm in $Z(\mathcal S^2)$ (see \cite[Proposition 7.2]{CCS} for a detailed proof or also the discussion in \cite{KM}). Therefore, we are clearly done. For the general case, given $n\in \mathbb N$, pick a subset $A\subseteq \mathbb N-\{2,...,2n\}$ with $|A|=n$ and then $\{v_j\}_{j\in A}$. Let $E=\{j: 2j\in A \}$. If $E=\emptyset$, we have that any element of $[v_j]_{j\in A}$ is of the form $(x,0)$ and the claim of the proposition is trivial. Otherwise, the set $E$ is admissible and thus, as before, $\Omega(y)=0$ if $\supp(y)\subseteq E$. If we pick $(x,y)\in [v_j]_{j\in A}$, it follows that $\supp(y)\subseteq E$ and thus $$\|(x,y)\|=\|x-\Omega(y)\|+\|y\|=\|x\|+\|y\|.$$ 
\end{proof}
However, $Z(\mathcal S^2)$ is not as. Hilbertian, so it is ``far" in some sense from $Z(\mathcal T^2)$ (see the picture in Section \ref{dibujito}). We need first a preparatory lemma.
\begin{lemma}\label{copia}
The twisted-Hilbert $Z(\mathcal S^2)$ contains an isomorphic copy of $Z_2$.
\end{lemma}
\begin{proof}
The blocks $u_n=2^{\frac{1-n}{2}}\sum_{j=2^{n-1}}^{2^n-1} e_j$, with $n\in \mathbb N$, span an isometric copy of $c_0$ in $\mathcal S^2$, see \cite[Proposition 0.7]{CS}. It is not hard to check that they span a copy of $\ell_1$ in $(\mathcal S^2)^*$. Indeed, since the support of each $u_n$ is an admissible set, we readily find that $\|u_n\|_{(\mathcal S_2)^*}=1$, and then
\begin{eqnarray*}
\left\| \sum_{n=1}^{N}\lambda_n u_n\right \|_{(\mathcal S^2)^*} &=& \left\| \sum_{n=1}^{N}|\lambda_n| u_n\right \|_{(\mathcal S^2)^*}\\
&=& \sup_{\|x\|_{\mathcal S^2}\leq 1} \sum_{n=1}^{N}|\lambda_n| \langle u_n, x\rangle\\
&\geq& \sum_{n=1}^{N}|\lambda_n|,
\end{eqnarray*}
where the last inequality holds picking the norm one vector $x=u_1+...+u_N$. Once this is achieved, it becomes trivial to check that given $(\lambda_n)_{n=1}^{\infty}\in c_{00}$ with $\|(\lambda_n)_{n=1}^{\infty}\|_{\ell_2}=1$ and letting $\lambda=\sum_{n=1}^{\infty}\lambda_n u_n$, the map
$$B(\lambda)(z):=\sum_{n=1}^{\infty}\lambda_n|\lambda_n|^{2z-1}u_n\in \mathcal F_{\infty}(\mathcal S^2, (\mathcal S^2)^*)$$
is well defined. It follows that $B(\lambda)$ is a $1$-bounded selection for $\lambda$ so that we find
$$\delta'_{1/2}B(\lambda)=\sum_{n=1}^{\infty}\lambda_n\log |\lambda_n|u_n.$$
The centralizer above induces the quasi-norm in the Kalton-Peck space (see \cite{KaPe}), thus $[(u_n,0),(0,u_n)]_{n=1}^{\infty}$ spans a copy of $Z_2$.
\end{proof}
%
We now prove our claim.
\begin{proposition}\label{z2}
The twisted-Hilbert $Z(\mathcal S^2)$ is not asymptotically Hilbertian.
\end{proposition}
\begin{proof}
First let us observe that $Z_2$ is not as. Hilbertian. Otherwise, the Orlicz space $\ell_M$ spanned by $(v_{2j})_{j=1}^{\infty}$ in $Z_2$ would be also as. Hilbertian: Obviously, if $X$ is finite codimensional in $Z_2$, then $X\cap \ell_M$ is also finite codimensional in $\ell_M$. But it is well known that as. Hilbertian spaces cannot have a symmetric basis unless they are isomorphic to $\ell_2$. Since $\ell_M$ is not isomorphic to $\ell_2$, otherwise $\ell_M$ would have type $2$ which is impossible by the unconditionality of the basis of $\ell_M$ since $$\left\|\left(0, \sum_{j=1}^n e_j \right) \right\|=\sqrt{n}\log{\sqrt{n}},$$ we are done. In particular, no space containing a copy of $Z_2$ is as. Hilbertian so that Lemma \ref{copia} allows us to conclude.
\end{proof}
 We now show that $Z(\mathcal S^2)$ is neither isomorphic to a subspace of the Kalton-Peck space nor a subspace of the Enflo-Lindenstrauss-Pisier space $\ell_2(\mathcal E_n)$ ( \cite{ELP}). We need to prove first the following lemma on the structure of $Z(\mathcal S^2)$.
\begin{lemma}\label{Pi3} There is $K>0$ such that for every finite subset $B\subseteq \mathbb N$, there is a subset $A\subseteq B$ with $|A|\geq 2^{-1}|B|$ for which $(v_j)_{j\in A}$ is $K$-equivalent to $\ell_2^{|A|}$.
\end{lemma}\begin{proof} Let $E=\{n:2n\in B\}$  and assume first that $|E|$ is even. Then, observe that the last $2^{-1}|E|$ elements of $E$ form an admissible set, namely $E'$. Thus for any $y\in c_{00}$ with $\supp (y)\subseteq E'$ we have
\begin{equation}\label{cdelta}
\max(\|y\|_{\mathcal S^2},\|y\|_{(\mathcal S^2)^*})=\|y\|_{\ell_2}.
\end{equation}
Let us show that the set $A:=\{2n-1\in B\}\cup \{2n: n\in E'\}$ satisfies the claim of the lemma. Given $y\in c_{00}$ as above we find as in Proposition \ref{Enn} $$\Omega(y)=0,$$ and hence, given $x\in\ell_2$ and $y\in c_{00}$ with $\supp (y)\subseteq E'$, we have that
$$\left \|x-\Omega(y) \right \| + \|y \|=\left \|x\right \|+\left \|y \right \|.$$
The expression above shows that $(v_j)_{j\in A}$ is equivalent, with uniform constant, to the basis of $\ell_2^{|A|}$. To finish observe that
$$|A|=|\{2n-1\in B \}|+ |\{2n: n\in E'\}|\geq |\{2n-1\in B \}|+ 2^{-1}|E|\geq 2^{-1}|B|,$$ and thus we are done. If $|E|$ is odd, then the last $(|E|+1)/2$ elements are again an admissible set and we argue as before.
\end{proof}
We now close the loop opened in Lemma \ref{copia}.
\begin{proposition} Let $\phi$ be a Lipschitz map such that either $\lim_{t\to \infty}\phi'(t)=0$ monotonically or $\phi(t)=ct$ for $c\neq 0$. Then, the following holds:
\begin{enumerate}
\item The space $Z(\mathcal S^2)$ is not isomorphic to either a subspace or a quotient of $Z_2(\phi)$.
\item The space $Z(\mathcal S^2)$ is not isomorphic to either a subspace or a quotient of $\ell_2(F_n)$ with $\dim F_n<\infty$ for all $n\in \mathbb N$.
\end{enumerate}
\end{proposition}
\begin{proof} The claims for the quotient maps follows by simple duality, for (1) see \cite[Lemma 1]{JS2}. The proof is similar to \cite[Proposition 4]{JS1} so we only sketch the argument. It is enough to show that no subsequence of $(v_{2j})_{j=1}^{\infty}$ in $Z(\mathcal S^2)$ is equivalent to $\ell_2$ or the Orlicz space $\ell_{M}$ of $Z_2(\phi)$ (see \cite[Lemma 5.3.]{KaPe}). If  $(v_{2n_j})_{j=1}^{\infty}$ is equivalent to $\ell_2$ then $(e_{2n_j})_{j=1}^{\infty}$ in $\mathcal S^2$ is also equivalent to $\ell_2$ by an easy application of Kalton's uniqueness theorem \cite[Theorem 7.6.]{ka}, and thus we trivially find a copy of $\ell_1 $ in $\mathcal S$ that is impossible \cite[Theorem 0.5]{CS}. Since the only symmetric basis satisfying Lemma \ref{Pi3} is $\ell_2$, the case of $\ell_{M}$ reduces to the previous one.
\end{proof}
The results of this section show that $Z(\mathcal S^2)$ is not isomorphic to $Z(\mathcal J)$ or to the previous examples of twisted Hilbert spaces: $\ell_2(\mathcal E_n)$, $Z_2$ and $Z(\mathcal T^2)$. Also, the results of Section \ref{ahnw} give that $Z(\mathcal J)$ is not isomorphic to $Z(\mathcal T^2)$. If, in addition, we use the fact that $Z_2$ is not as. Hilbertian (given in the proof of  Proposition \ref{z2}), then $Z(\mathcal J)$ is also not isomorphic to $Z_2$. It only remains to check that $Z(\mathcal J)$ is different from $\ell_2(\mathcal E_n)$. This follows from our next proposition which needs a preparatory result.
\begin{lemma}\label{critico} Let $X$ denotes a space with a shrinking FDD, say $(F_n)_{n=1}^{\infty}$, and let $F$ be an $m$-codimensional subspace of $X$. Given $\varepsilon>0$, there is $M=M(m,\varepsilon)\in \mathbb N$ so that $F$ contains an $(1+\varepsilon)$-isomorphic copy of $[F_n]_{n=M}^{\infty}$.
\end{lemma}
\begin{proof}
For a given $\delta>0$ (that will be fixed later), we may find in the unit sphere of the $m$-dimensional space $(X/F)^*$ a finite $\delta$-net, say $\{y_j^*\}_{j=1}^N$, just by compactness.

Now, for each $j\leq N$, by the shrinking condition on the $(F_n)_{n=1}^{\infty}$, we may pick $m_j\in \mathbb N$ so that $$\|y_j^*\circ Q_{|[F_n]_{n=m_j}^{\infty}}\|\leq \frac{\varepsilon}{4m},$$ 
where $Q:X\longrightarrow X/F$ stands for the natural quotient map. Picking $M=\max_{j\leq N} m_j$, let us check that 
\begin{equation}\label{ringurango}
\|Q(x)\|\leq \frac{\varepsilon}{2m}\|x\|, \;\;x\in [F_n]_{n=M}^{\infty}.
\end{equation}
Indeed, pick $x\in [F_n]_{n=M}^{\infty}$ with $\|x\|=1$. Then $\|Qx\|\leq 1$ and we may choose $x^*\in (X/F)^*$ also with $\|x^*\|=1$ so that $\|Q(x)\|=|x^*(Q(x))|$. For the choice $x^*$ there is $j_0\leq N$ so that $\|x^*-y_{j_0}^*\|\leq \delta$ by the very definition of $\delta$-net. Thus,
\begin{eqnarray*}
\|Q(x)\|&=&|x^*(Q(x))|\\
&\leq&|(x^*-y_{j_0}^*)(Q(x))|+|y_{j_0}^*(Q(x))|\\
&\leq& \delta+\frac{\varepsilon}{4m}.
\end{eqnarray*}
If we choose $\delta=\frac{\varepsilon}{4m}$, we may take for granted that (\ref{ringurango}) is achieved. 

We are going to use the estimate $(\ref{ringurango})$ for the bound of a projection of $X$ onto $F$. The construction of the projection is standard. Pick an Auerbach basis in $X/F$, say $\{z_1,...,z_m\}$ with biorthogonal functionals $\{z_1^*,...,z_m^*\}$ so that $z_i^*(z_j)=\delta_{ij}$ and $\|z_j\|=\|z_j^*\|=1$ for $j\leq m$. For each $z_j$ pick $x_j\in X$ so that $Q(x_j)=z_j$ and $\|x_j\|\leq 2$ for all $j\leq m$. Therefore, the map
$$P:X\longrightarrow \Ker Q=F,$$
given by the rule $$P(x)=x-\sum_{j=1}^m\langle z_j^*, Q(x) \rangle x_j$$
is well defined since $$Q\left( \sum_{j=1}^m\langle z_j^*, Q(x) \rangle x_j \right)=\sum_{j=1}^m\langle z_j^*, Q(x) \rangle z_j=Q(x).$$
And it is bounded since
\begin{eqnarray*}
\|P(x)\|&=&\left\|x-\sum_{j=1}^m\langle z_j^*, Q(x) \rangle x_j\right\|\\
&\leq& \|x\|+\sum_{j=1}^m|\langle z_j^*, Q(x) \rangle|\| x_j\|\\
&\leq&\|x\|+\|Q(x)\|2m\\
&\leq& (1+2m)\|x\|.
\end{eqnarray*}
But let us observe that the estimate $(\ref{ringurango})$ gives for $x\in [F_n]_{n=M}^{\infty}$ that $$\|P(x)\|\leq (1+\varepsilon)\|x\|,$$
and also $$\|P(x)\|\geq \|x\|-2m\|Q(x)\|\geq (1-\varepsilon)\|x\|.$$
The previous two estimates show that $P$, when restricted to $[F_n]_{n=M}^{\infty}$, is an almost isometry provided $\varepsilon>0$ is small enough.
\end{proof}
A similar claim for spaces with a shrinking basis seems to be folklore but we found no proof in the literature so we have included one for the sake of completeness. Now, we may conclude our previous discussion and distinguish between $\ell_2(\mathcal E_n)$ and the asymptotically Hilbertian $Z(\mathcal J)$.
\begin{proposition}
The space $\ell_2(\mathcal E_n)$ is not asymptotically Hilbertian.
\end{proposition}
\begin{proof}
Let us recall a couple of useful facts on $\ell_2(\mathcal E_n)$. The first is that $\mathcal E_n=\ell_2^{n^2}\oplus_{f_n}\ell_2^n$ where $f_n$ is a certain nonlinear function. The key is that such function is defined recursively, so that given $f_n$ we construct $f_{2n}$. Using this fact and the expression \cite[(12), Lemma 1]{ELP}, one trivially finds that 
\begin{equation}\label{elpf}
\mathcal E_n \subseteq \mathcal E_{2n}
\end{equation}
isometrically for every $n\in \mathbb N$. A second useful fact is that 
\begin{equation}\label{estimation}
d_{\mathcal E_n}\geq c\sqrt {\log n},
\end{equation}
for some absolute $c>0$. This follows trivially from \cite[Corollary]{ELP}. Let us sketch now that $\ell_2(\mathcal E_n)$ is not as. Hilbertian. Let us fix the dimension $n^2+n$ and pick a finite codimensional subspace $F$ of $\ell_2(\mathcal E_n)$ so that every $(n^2+n)$-dimensional subspace of $F$ is $C$-isomorphic to a Hilbert. By Lemma \ref{critico}, since $\ell_2(\mathcal E_n)$ has clearly a shrinking FDD given by $(\mathcal E_n)_{n=1}^{\infty}$, such $F$ contains a $2$-isomorphic tail of the FDD, say $(\oplus_{j=g(n)}^{\infty}\mathcal E_j)_2$. So the same claim must also hold for this tail of the FDD but with constant $2C$. It is clear by (\ref{elpf}) that one may find $m\geq g(n)$ large enough so that $$\mathcal E_n\subseteq \mathcal E_m,$$
where let us recall that $\dim \mathcal E_n=n^2+n$. But using (\ref{estimation}) we find that $c\sqrt {\log n}\leq d_{\mathcal E_n}\leq 2C$ which is absurd if $n$ is large enough.
\end{proof}
\section{A new HAPpy space}\label{cinco}
We construct now a nontrivial twisted-Hilbert which is a HAPpy space. The main technical tool to show it is HAPpy is a result of Johnson and Szankowski claiming that if $d_n(X)$ grows to infinity slowly enough, of inverse Ackermann type, then $X$ is HAPpy, see \cite[Theorem 2.1]{JS}. So let us introduce the inverse Ackermann function. We must first familiarize ourselves with the hierarchy of rapidly growing functions. We define 
$$\left \{\begin{array}{l l}g_0(k)=k+1, & \mbox{and } \\g_{n+1}(k)=g^{(k)}_n(k), & \mbox{for }n\geq 0,\end{array}\right. $$
where $g^{(k)}_n$ denotes the $k$-fold iterate of $g_n$, i.e., $g^{(i+1)}_n(k)=g_n(g^{(i)}_n(k))$. It follows that 
$$\left \{\begin{array}{l l}g_1(k)=2k, & \\g_2(k)=k2^k, & \end{array}\right. $$
while the map $g_3(k)$ corresponds to a stacked tower of height $k$. The inverse of the Ackermann function is the inverse of the map $$n\mapsto g_n(n).$$
For simplicity we shall deal with the inverse of the map $$n\mapsto g_n(2),$$
which, as was pointed in \cite{JN}, has the same asymptotic behaviour that of the inverse Ackermann function. This is, we define this inverse $\alpha(n)$ to be the unique integer $i$ so that $$g_{i}(2)\leq n< g_{i+1}(2).$$ The reader will find that any other variant of the inverse Ackermann function is the same up to a bounded factor \cite[Lemma B.1., Appendix B]{AKa}. As an example of the appearing of the hierarchy of rapidly growing functions in the literature we have the following well known result \cite{CS}.
\begin{proposition}
There is $C>0$ so that every $g_n(k)$-dimensional subspace of $[t_j]_{j=k}^{\infty}$ in $\mathcal T^2$ is $C^n$-isomorphic to a Hilbert space.
\end{proposition}
In particular, we may describe the Banach-Mazur distance to a Hilbert space in terms of the inverse Ackermann function.
\begin{corollary}
For $n\in \mathbb N$ large enough, every $g_{\alpha(n)+1}(2)$-dimensional (thus every $n$-dimensional) subspace of $\mathcal T^2$ is $2^{O(\alpha(n))}$-isomorphic to a Hilbert space.
\end{corollary}
In this section we deal with the Casazza-Nielsen symmetric version of $\mathcal T^2$ that we will denote by $\mathcal T^2_s$. The reader may fulfill all the details from the paper of Casazza and Nielsen \cite{CN2} where this space is denoted by $S(\mathcal T^2)$. To simplify things, we just use two facts of this space. One aspect is that the natural basis of $\mathcal T^2_s$ is symmetric. The second ingredient is that Casazza and Nielsen proved that $d_n(\mathcal T^2_s)\to \infty$ more slowly than any iteration of the logarithm \cite[Proposition 3.9]{CN2}. The proof indeed shows that $d_n(\mathcal T^2_s)$ is of inverse Ackermann type as in the corollary above.
\begin{proposition}\label{Casazza} For every finite dimensional subspace $E$ of $\mathcal T^2_s$ of large enough dimension,
\begin{equation}\label{14}
d_E\leq 2^{O(\alpha(\dim E))}.
\end{equation}
\end{proposition}
\begin{proof} We only sketch it since the argument is the same as in \cite[Proposition 3.9]{CN2}. Let $n=\dim E$ and using a standard argument we may assume that we are on the span of $n^n$ vectors with disjoint support. Then, arguing exactly as in \cite[Proposition 3.9]{CN2}, and using the corollary above with estimate $2^{c\alpha(n)}$ for $n$-dimensional subspaces where $c>0$ is some fixed absolute constant, we infer that $$d_E\leq K\cdot 2^{2c\alpha(n^n)},$$
for some absolute $K>0$.
We only need to check that $$\alpha(n^n)\leq \alpha(n)+2,$$
but this follows from
\begin{eqnarray*}
g_{\alpha(n^n)}(2)\leq n^n&\leq &g_3(n)\\
&\leq& g_3(g_{\alpha(n)+1}(2))\\
&\leq& g_{\alpha(n)+1}(g_{\alpha(n)+1}(2))\\
&=&g_{\alpha(n)+2}(2),
\end{eqnarray*}
where we have assumed that $2\leq \alpha(n)$ in the fourth inequality.
\end{proof}
 A similar estimate holds for $(\mathcal T^2_s)^*$. We may give a short argument (that will be useful later) based on an idea of Johnson \cite{BJ}.
\begin{corollary}\label{corcasazza}
For every finite dimensional subspace $E$ of $(\mathcal T^2_s)^*$ of large enough dimension,
\begin{equation}\label{15}
d_E\leq 2^{O(\alpha(\dim E))}.
\end{equation}
\end{corollary}
\begin{proof}
Fix such $E$ and pick $F$ in $\mathcal T^2_s$ that is $2^{-1}$-norming. This can be achieved easily with $\dim F\leq 5^{\dim E}$, see for example \cite[Lemma, page 7]{SM}. Then the natural duality pairing $B:E\times F\to \mathbb K$ given by $B(e,f)=e(f)$ provides a $2$-isomorphic embedding of $E$ into $F^*$. Therefore,
$$d_E\leq 2d_{F^*}=2d_F\leq 2\cdot2^{c\alpha(5^{\dim E})}\leq 2^{c_1\alpha(\dim E)},$$
where in the last inequality we have argued as in the proof of Proposition \ref{Casazza} for $\dim E$ large enough and in the previous inequality we have used the estimate guaranteed by Proposition \ref{Casazza} for some $c>0$.
\end{proof}
The $\ell_2$ basis dominates the basis in $\mathcal T_s^2$, therefore $(\mathcal T_s^2, (\mathcal T_s^2)^*)_{1/2}=\ell_2$ by \cite[Corollary 4.3]{CoS}. Let $Z(\mathcal T_s^2)$ be,  as usual, the induced twisted Hilbert space. Thus, we are ready to prove the twisted analogue of the result of Casazza-Nielsen.
\begin{proposition}\label{distancia}For every finite dimensional subspace $E$ of $Z(\mathcal T^2_s)$ of large enough dimension,
\begin{equation}\label{16}
d_E\leq 2^{O(\alpha(\dim E))}.
\end{equation}
\end{proposition}
\begin{proof}
Let us observe the trivial bound $a_{n,2}(X)\leq d_n(X)$. Therefore, using the estimates (\ref{14}) and (\ref{15}), we have, for some absolute $c>0$, $$\max\{a_{n,2}(\mathcal T^2_s), a_{n,2}((\mathcal T^2_s)^*)\}\leq 2^{c\alpha(n)},$$
if $n$ is large enough. We use now the estimate provided in \cite[Proposition 3]{JS1} to find
$$a_{n,2}(Z(\mathcal T_s^2))\leq c_1\cdot 2^{c\alpha(n)},$$
for some absolute $c_1>0$. If we let now $E$ to be an $n$-dimensional subspace, then, reasoning exactly as in Claim B of Proposition \ref{as}, we find:
\begin{eqnarray*}
d_E&\leq& a_2(E)\cdot c_2(E)\\
&\leq&2\sqrt{2\pi}\cdot a_{n,2}\left(Z(\mathcal T_s^2) \right)\cdot c_{n,2}\left( Z(\mathcal T_s^2) \right)\\
&\leq& 2\sqrt{2\pi}\cdot \lambda \cdot a_{n,2}\left(Z(\mathcal T_s^2)\right)^2\\
&\leq& 2\sqrt{2\pi}\cdot \lambda \cdot c_1^{2}2^{2c\alpha(n)},  
\end{eqnarray*}
where $\lambda$ denotes the isomorphism constant between $Z(\mathcal T_s^2)$ and its dual \cite[Lemma 1]{JS2}.
\end{proof}
\begin{corollary}
The spaces $Z_2, \ell_2(\mathcal E_n)$ and $Z(\mathcal S^2)$ are not isomorphic to a subspace or a quotient of $Z(\mathcal T^2_s)$.
\end{corollary}
To state our next corollary let us recall that Johnson and Szankowski proved that if $d_n(X)\to \infty$ slowly enough (of inverse Ackermann type) then $X$ is HAPpy, see \cite[Theorem 2.1]{JS}. Therefore, it readily follows from the estimate (\ref{16}) that
\begin{corollary}
The twisted-Hilbert $Z(\mathcal T_s^2)$ is a HAPpy space.
\end{corollary}
Previously, the weak Hilbert $Z(\mathcal T^2)$ was the only nontrivial example of a twisted Hilbert that is HAPpy. It only remains to separate the new guy $Z(\mathcal T_s^2)$ from the as. Hilbertian gang.
\begin{proposition}\label{propsym}
The space $Z(\mathcal T^2_s)$ is not asymptotically Hilbertian nor the basis $(v_j)_{j=1}^{\infty}$ has the $E(n,n,K)$-property.
\end{proposition}
\begin{proof}
Let us assume that $Z(\mathcal T^2_s)$ is as. Hilbertian and let us reach a contradiction. Since the basis of $\mathcal T_s^2$ and its dual are symmetric, there is no loss of generality assuming that the corresponding centralizer $\Omega$ is symmetric. This means exactly that there is $C>0$ such that for every permutation $\sigma$ of the natural numbers, we have that 
\begin{equation}\label{symmkalt}
\|T_\sigma \circ \Omega(x)-\Omega\circ T_\sigma (x)\|\leq C\|x\|,\;\;\;x\in c_{00},
\end{equation}
where $T_{\sigma}(\sum_{n=1}^{\infty} a_n e_n)=\sum_{n=1}^{\infty} a_n e_{\sigma(n)}$. Indeed, pick $B$ a $2$-bounded selector so that $\delta_{1/2}'B=\Omega$ and observe that the linearity of $T_{\sigma}$ easily implies  $$T_{\sigma}\circ \delta_{1/2}'B=\delta_{1/2}'T_{\sigma}\circ B.$$ Therefore, 
$$T_\sigma \circ \Omega-\Omega \circ T_\sigma =\delta_{1/2}'(T_{\sigma}\circ B-B \circ T_{\sigma}).$$
The key point is that $T_{\sigma}\circ B(x)-B \circ T_{\sigma}(x)\in \Ker \delta_{1/2}$ so that we may apply the compatibility criteria cf. \cite[Theorem 4.1.]{CCS}, that is
$$\|\delta_{1/2}'(T_{\sigma}\circ B(x)-B \circ T_{\sigma}(x))\|_2\leq C \| T_{\sigma}\circ B(x)-B\circ T_{\sigma}(x)\|_{\mathcal F_{\infty}},$$
for some absolute $C>0$. It is trivial to get a good bound for $B\circ T_{\sigma}$. Now, to bound $T_{\sigma}\circ B$ we must use that the operators $T_{\sigma}$ are (uniformly) bounded when acting on $\mathcal T_s^2$ and its dual since the basis are symmetric. We leave the easy details to the reader. Just to mention, Kalton proved the very much difficult converse to this claim, namely, if the centralizer is symmetric then we may assume that the spaces representing such centralizer have a symmetric basis \cite[Corollary 7.7]{ka}.

It follows trivially from (\ref{symmkalt}) that the basis $(v_{2j})_{j=1}^{\infty}$ is symmetric. The only as. Hilbertian with a symmetric basis is $\ell_2$ but if $[v_{2j}]_{j=1}^{\infty}\approx \ell_2$ then, once more, by Kalton's uniqueness theorem \cite[Theorem 7.6.]{ka}, we would find $\mathcal T_s^2\approx \ell_2$. A very similar argument holds for the $E(n,n,K)$-property.
\end{proof}
\begin{corollary} The following holds:
\begin{enumerate}
\item $Z(\mathcal T^2_s)$ is not isomorphic to a subspace or a quotient of $Z(\mathcal T^2)$.
\item $Z(\mathcal T^2_s)$ is not isomorphic to a subspace or a quotient of $Z(\mathcal J)$.
\end{enumerate}
\end{corollary}
\subsection{The Johnson-Szankowski twisted Hilbert spaces}
This subsection deals with another theorem of Johnson and Szankowski \cite[Theorem 3.1.]{JS}:
\begin{theorem}\label{joh}
Let $1<\delta_n\to \infty$. There exists an Orlicz  space $\ell_M$ of type 2, non isomorphic to $\ell_2$, so that $d_n(\ell_M)\leq \delta_n$ for every $n\in \mathbb N$.
\end{theorem}
The arguments given in Section \ref{cinco} can be used to prove a ``twisted analogue" of this result. In particular, it is a way to produce HAPpy twisted Hilbert spaces.
\begin{proposition}\label{twjoh} Let $1<\delta_n\to \infty$. There exists $n_0\in \mathbb N$ and a twisted Hilbert space $Z(\mathcal {JS})$, so that $d_n(Z(\mathcal {JS}))\leq \delta_n$ for $n\geq n_0$. Moreover, $Z(\mathcal {JS})$ is not asymptotically Hilbertian.
\end{proposition}
\begin{proof} The proof uses the same arguments given for $Z(\mathcal T_s^2)$. Let $(\delta_n)_{n=1}^{\infty}$ be fixed and pick $(\varepsilon_n)_{n=1}^{\infty}$ with $1<\varepsilon_n\to \infty$ so that $$\eta \cdot (2\cdot \varepsilon_{5^n})^2\leq \delta_n,$$ where $\eta$ is a universal constant (to be discussed later) and $n$ is large enough: For example, define ``$\varepsilon_n$" to be the constant function $\sqrt{\delta_n}/(4\eta)$ on each exponential jump. Consider now the Johnson-Szankowski space of Theorem \ref{joh} for the choice $(\varepsilon_n)_{n=1}^{\infty}$, say $\ell_M$. Then, let us pick $Z(\mathcal {JS}):=Z(\ell_M)$. By construction $$d_n(\ell_M)\leq \varepsilon_n$$ and, using the argument of Corollary \ref{corcasazza}, we also have $$d_n(\ell_M^*)\leq2\cdot d_{5^n}(\ell_M)\leq  2\cdot \varepsilon_{5^n}.$$ Let us argue as in Proposition \ref{distancia} to find 
\begin{equation}\label{last}
d_n(Z(\ell_M))\leq \eta\cdot (2\cdot \varepsilon_{5^n})^2,
\end{equation}
 for some constant $\eta$ which deserves a couple of remarks. It appears as a consequence of two facts:
\begin{itemize}
\item The duality: the isomorphism constant of $Z(X)$ with its dual does not depend on $X$. A detailed proof is to be found in \cite[Lemma 1]{JS2} where it is shown that it only depends on the constant appearing in (\ref{centralizer}) which in turn is bounded by $4\rho$; recall that $\Omega=\delta_{1/2}'B$ where $B$ is a $\rho$-bounded selector. Since we may clearly assume without loss of generality that $\rho=2$, the isomorphism constant can be bounded independently of $X$.
\item The type $2$ constants: A close inspection of \cite[Proposition 3]{JS1} shows that if we denote by $$A_n=\max\{a_{n,2}(X), a_{n,2}(X^*)\},$$ then we have that $$a_{n,2}(Z(X))\leq 4\rho \cdot A_n.$$
And, as before, we may assume that $\rho=2$.
\end{itemize}
 Therefore, it follows from the comments above that there is such $\eta$ (independent of the space involved) so that (\ref{last}) holds. Indeed, repeating the chain of estimates provided for $d_E$ in Proposition \ref{distancia} with this new, and much more precise, labelling of the constants, we find for an $n$-dimensional subspace $E$ of $Z(X)$ that $$d_E\leq \eta \cdot A_n^2,$$ where $\eta$ is independent of $X$. Thus, we have $d_n(Z(X))\leq\eta\cdot A_n^2$ and then also
$$
d_n(Z(X))\leq \eta \cdot D_n^2,
$$
where $$D_n=\max\{d_n(X), d_n(X^*)\}.$$ Finally, using the technique of Corollary \ref{corcasazza}, we arrive to the easy-to-handle formula $$d_n(Z(X))\leq 4\eta\cdot 5^{2d_n(X)}.$$ Hence, we are done with the first part. To finish we need to check that $Z(\ell_M)$ is not as. Hilbertian but this is exactly as in Proposition \ref{propsym}. 
\end{proof}
\begin{corollary}
Let $1<\delta_n\to \infty$. There exists $n_0\in \mathbb N$ and a $\ell_2$-saturated space $X$ with a symmetric basis but failing type $2$ so that $d_n(X)\leq \delta_n$ for $n\geq n_0$.
\end{corollary}
\begin{proof}
Pick the twisted Hilbert $Z$ of Proposition \ref{twjoh} for the choice $(\delta_n)_{n=1}^{\infty}$ and consider as $X$ the closed span of $(v_{2j})_{j=1}^{\infty}$.
\end{proof}
\section{A pocket map for the hexagonal tree}\label{dibujito}
The following picture may help to understand the relationship of the $6$-known types of  twisted Hilbert spaces so far.  

The basic idea is that every space is connected to another if they are similar in some sense. For example, $Z_2(\phi)$ is connected to $\ell_2(\mathcal E_n)$ since the distance of its $n$-dimensional subspaces to Hilbert is extremal, roughly, of the order $\log n$.  Another one is that $Z(\mathcal J)$ and $\ell_2(\mathcal E_n)$ are connected since both spaces are $\ell_2$-sums of finite dimensional spaces. The link of $Z(\mathcal J)$ and $Z(\mathcal T^2)$ is that both of them are asymptotically Hilbertian while this last and the Hilbert copy are both weak Hilbert spaces. On the other hand, $Z(\mathcal T^2)$ shares with $Z(\mathcal T^2_s)$ that both are HAPpy while it shares with $Z(\mathcal S^2)$ that both do satisfy the $E(n,n,K)$-property. To close the loop, $Z(\mathcal S^2)$ is connected to $Z_2$ since this last is a subspace of $Z(\mathcal S^2)$. In a more elaborated language, this connection may be stated as both twisted Hilbert spaces are admitting $Z_2$ as a ``twisted" spreading model, something that does not hold for the rest of the spaces of the tree. These claims on spreading models may be found in the forthcoming paper \cite{JS4}.

We have implemented the intuitive idea of the ``size" of such similarity in terms of the distance. In this sense, $\ell_2$ plays the role of the tree root, it is very close to $Z(\mathcal T^2)$ (both of them are weak Hilbert spaces) but as far as possible from $Z_2$ which is, undoubtedly, the treetop.
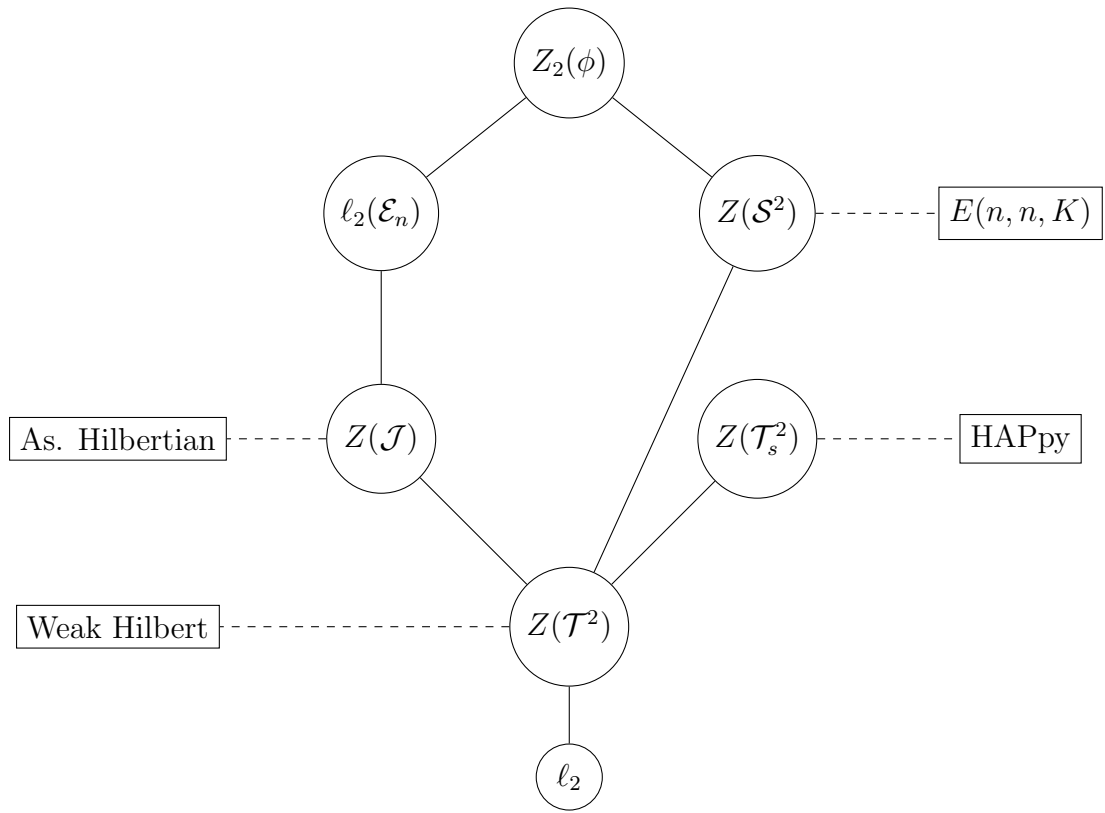
\begin{center}
\begin{tikzpicture}[square/.style={regular polygon,regular polygon sides=4}]
\path (0,0) node[draw,shape=circle] (p0) {$Z_2(\phi)$}
(-2.5,-2) node[style={draw,shape=circle}] (p1) {$\ell_2(\mathcal E_n)$}
(-2.5,-5) node (p2)[style={draw,shape=circle}] {$Z(\mathcal J)$}
(2.5,-5) node (p3)[style={draw,shape=circle}] {$Z(\mathcal T_s^2)$}
(2.5,-2) node (p4)[style={draw,shape=circle}] {$Z(\mathcal S^2)$}
(0,-7.5) node (p5)[style={draw,shape=circle}] {$Z(\mathcal T^2)$ }
(0,-9.5) node (p6)[style={draw,shape=circle}] {$\ell_2$}

(-6,-7.5) node (p9)[style={draw,shape=rectangle}] {Weak Hilbert}
(-6,-5) node (p10)[style={draw,shape=rectangle}] {As. Hilbertian}

(6,-2) node (p11)[style={draw,shape=rectangle}] {$E(n,n,K)$}
(6,-5) node (p12)[style={draw,shape=rectangle}] {HAPpy};
\draw 
(p0) -- (p1)
(p0) -- (p4)
(p1) -- (p2)
(p2) -- (p5)
(p3) -- (p5)
(p4) -- (p5)
(p5) -- (p6);
\draw[dashed]
(p2) -- (p10)   
(p5) -- (p9)    
(p4) -- (p11)   
(p3) -- (p12);
\end{tikzpicture}
\end{center}

\begin{remark} All our examples come from the complex interpolation method. The spaces $Z_2, Z(\mathcal T^2), Z(\mathcal J)$ and $Z(\mathcal T_s^2)$ admit a version which is given by the real interpolation method, so the constructions are not attached to one and only method of interpolation. The interested reader may find this in \cite{JS3}. 
\end{remark}
\textbf{Acknowledgement.}
We are extremely grateful to the referee for the careful reading of the manuscript. It leads us to a substantial improvement in the quality and the presentation of the paper.

\end{document}